# The family of Eulerian constants in the isodimensional discrete calculus

Ricardo Ramos Montero

*Departamento de Informática, Universidad de Oviedo*
Campus de Viesques, 33204 Gijón, Asturias, Spain

rramos@uniovi.es

Discrete Euclidian Spaces (DESs) are the starting point in the study of the major fields of the Isodimensional Discrete Mathematics (IDM). The *isodimensional analysis* is not an exception, being particularly interesting and fruitful the presence of the DESs in the study of differential and integral calculus. At the moment, and because of the novelty, the discrepancies between the *standard* and the *isodimensional discrete calculus* (IDC) are what draw more attention. As a small sample, in this paper we will see the interpretation given by the IDC about the origin and nature of the constant of Euler, and after that, we will also do likewise with another pair of constants of the same family.

## Introduction

The *isodimensional discrete calculus* (IDC)— which belongs to the *isodimensional discrete mathematics* (IDM), and whose analysis is the primary objective of this paper —is practically a newcomer, as it cannot be found as such in traditional mathematics, and however, it can be said that the IDC hardly has its own concepts. Then, what are its merits?

The main contribution that the IDC makes to the *discrete analysis* is a common framework, both valid for the discrete version of the *infinitesimal calculus*, and for the set of techniques focused to work with finite elements, such as the *finite differences methods*. This grouping of methodologies in a common framework brings the unification of the basic ideas and terminology. In short, the great merit of the IDC is to provide an exclusive infrastructure, which allows the compilation and unification of concepts and techniques of calculus, providing in this way a more intuitive overview of the differential and integral calculus. Is there something similar in mathematics?

The answer is affirmative, since there are at least two areas of mathematics that could be confused with the IDC. On the one hand we have the *discrete calculus* [1][2] (also known as *discrete exterior calculus)*, which has its roots in *graph theory*. It was born more or less at the same time as the *vector calculus*, so it cannot be said to be new. A prominent feature of the discrete exterior calculus is that it is fully referenced on a discrete space. This means that its goal is not to discretize the *standard* [infinitesimal] *calculus*, but to establish a genuine discrete one. As the IDC takes its main ideas and concepts from the infinitesimal calculus and also from the finite differences methods, the conceptual differences we can find between them are those that there are between the standard calculus (and the finite differences methods), and the discrete exterior calculus. In any case, they have in common that both work in their respective *discrete spaces,* i.e., neither of them try to discretize the infinitesimal calculus.

On the other hand there is the *time scale calculus*, which took its first steps in 1988, in the hands of Stefan Hilger [3]. One of the primary objectives of this calculus is the unification of the differential and integral calculus [4], with the calculus of differences, i.e., it has more or less the same goals as the IDC, only that the *time scale calculus* uses for its purposes the temporal scale (as its name suggests), while the IDC works with the spatial one[1].

To appreciate the IDC in detail, it is necessary to know the theoretical context where it develops[2], that is, the fundamentals and basic concepts of the IDM, although a summary of them [5] is sufficient to understand some key aspects of the isodimensional calculus. Taken for granted this minimum background, in [5] it can be seen that at least two scalar levels are involved in the definition of numbers. On one side, it is the *level of evaluation*, which defines or sets the value of numbers, and on the other, it is the *level of observation* which, as its name suggests, is the level from which we consider the value of a number. For example, in the *evaluated numerical sequence* 034.2871,3..., the comma indicates the level of evaluation *($S_6$)*, and the decimal point shows the perception level *($S_2$)*.

---

[1] So, another possible name for the IDC would be the *space scale calculus*.

[2] All material developed on the IDM is collected in a mathematical treatise, not yet published (www.isodimensional.com).

Then, extrapolating these ideas to the differential and integral calculus, the *level of perception* will usually be indicated by $S_0$, and the *level of definition* by $S_m$ (with $m \geq 0$), i.e., $S_m$ will be the scalar level (specific or generic) at which the calculation takes place. With these simple premises we are able to understand some basic issues of the IDC. Let us first see how (and with what ease) the IDC integrates in the same context the standard calculation and other mathematical techniques, such as the finite difference methods.

## Local and extreme calculus

Whatever the type or nature of the calculation we make, it must be clear that it is always carried out at the $S_m$ level. In addition, at this level, as in any other, the points of the local DES are discrete, which implies that the minimal variation or increment of a variable is 1.

Without losing sight of this pair of key ideas, the *quotient of increments* of the function $y = f(x)$, which leads directly to the concept of derivative function is given in the *isodimensional discrete analysis* by the generic algebraic expression[1]

$$\frac{\Delta y}{\Delta x} = \frac{f(x+\Delta x)-f(x)}{(x+\Delta x)-x} = \frac{f(x+\Delta x)-f(x)}{\Delta x} = f'(x) \qquad \text{E. 1}$$

where $\Delta x = uh$, with $h = 1/b^m$, and $u, b, m \in N$, being $b \geq 2$ the scalar order, and $u \geq 1$. ¿Why $\Delta x = uh$?

First of all, $u$ is an arbitrary integer to give generality to the expression. However, typically we are interested in a minimal variation of $\Delta x$, therefore $u = 1$ in most cases. So, assuming from now on that $u = 1$, we have that $\Delta x = h$, being $h = 1/b^m$ the *minimal subscalar variation*[2], whose value, as we see, depends on $m$ (the *scalar gap* or *height*), i.e., the number of scalar levels that there are from the evaluation level $(S_0)$, down to the level of definition of the calculation $(S_m)$.

According to the above, if the scalar gap is zero $(m = 0)$, i.e., if we look at the calculation on the same level at which it runs $(S_m = S_0)$, then $\Delta x = 1/b^0 = 1$, which corresponds to the minimum increment in any local DES, as expected.

---

[1] Note that the concept of limit is not involved at all in the generic definition of the derivative, due to the discrete nature of the IDM.

[2] The *minimal subscalar variation* plays a key role throughout the IDM, not only in the isodimensional calculus.

This means that if the calculation is observed from its own level of definition, then the derivatives (forward and backward) are given by $f_+(x) = f(x + 1) - f(x)$ and $f_-(x) = f(x) - f(x - 1)$, respectively[1]. On the other hand, if we assume that the scalar gap grows limitlessly $(m \to \infty)$[2], then $\Delta x \to 0$. We then see that if the $S_0$ level is too far from $S_m$, the derivatives and integrals in the IDC are in practice equal to those found in the standard calculus, although their definition is conceptually different. In other words, with $m = 0$ we have a local vision of the differential and integral calculus, while doing $m \to \infty$, we have a global one, although the calculation is carried out at $S_m$ in both cases. What if $m = k$, being $0 < k << \infty$? In this case we would have an *intermediate* type of calculus, which must be developed for specific problems or situations, depending upon the value of $k$.

In summary, the DESs provide a set of "systems of differential and integral calculus" in the scalar range $[S_0, S_m]$, with $m \geq 0$. From a mathematical point of view, the most important calculus systems are without doubt, those with $m = 0$ or $m \to \infty$, because they are generic. According to the IDC terminology, the *extreme* or *global calculus*— which in practice is equal to the standard calculus —has the perception level $S_0$ far from $S_m$ $(m \to \infty)$, and the *local calculus*— which embraces all those discrete techniques where it happens that $\Delta x \geq 1$ — has it at $S_m$ $(m = 0)$. In any case, the calculation always takes place at $S_m$.

Note that it would not make much sense to change the concepts involved, or the calculus terminology, because of varying the level of observation $(S_0)$. Hence, the greatest virtue of the IDC is the unification of techniques, ideas and terms, that are equal in all systems of calculus in the range $[S_0, S_m]$, except for a small change in the terminology, that somehow must indicate the scalar gap chosen.

## Local derivatives of the logarithmic function

We have just seen the general ideas that characterize the IDC and so now is the time to implement them. To do this, nothing better than analyzing a particular case, such as the logarithmic function.

---

[1] Similar expressions can be found in the finite differences methods.

[2] In the IDM the symbol '$\infty$' represents the *discrete infinite*, a concept similar to that used in mathematics when it defines the sequence of natural numbers.

Let $F(x) = \ln x$. Assuming that $m = 0$, i.e., making the evaluation of the calculation at the same scalar level where it is carried out, we will search for the *forward derivative* of $F(x)$, which is given by $f_+(x) = \ln (x + 1) - \ln x = \ln ((x + 1)/x)$, as we saw above. This is the *exact local f-derivative*, so named to differentiate it from the *approximate local f-derivative*, which will be indicated by $f_+(x)_\sim$. How is this derivative defined?

Being $f_+(x) = \ln ((x + 1)/x) = u$, it happens by definition that $e^u = (x + 1)/x$, and so, $e^u = (1 + 1/x)$. Raising both sides to the integer power $x$, we have that $e^{ux} = (1 + 1/x)^x$. If we now assume that $x \to \infty$, then $e^{ux} = (1 + 1/x)^x = e$. I must clarify that in the IDM, to calculate $e$ in this way, $x$ must take a value as large as we want (depending on the precision that we want to give to $e$), but in any case $x$ must be *finite*[1]. Therefore, we assume that $t$ is an integer value of $x$, from which $(1 + 1/x)^x$ provides the value of $e$ with the required precision. Consequently, $(1 + 1/x)^x = e \Leftrightarrow x \geq t \gg 1$, and $(1 + 1/x)^x < e$, if $x < t$.

Since $e^{ux} = e$, necessarily $ux = 1$, so that $u = 1/x$. Taking this result above finally we have that

$$f_+(x)_\sim = 1/x \qquad \text{E. 2}$$

Why do we say that this *f-derivative* is approximate? Simply because the numerical results that it provides are only reliable in the interval $[t, \infty]$. For $x < t$, the calculation is imprecise, since the algebraic development that leads to **E. 2** is incorrect, seeing that $e^{ux} < e \Rightarrow u \neq 1/x$. In short, and to be a little more concise, $f_+(x) = \ln ((x + 1)/x)$ is the *exact local f-derivative* in the range $[1, t]$, while $f_+(x)_\sim = 1/x$ is the *approximate local f-derivative* in the same interval.

Since the *extreme* (or *global*) *f-derivative* of $F(x) = \ln x$ is also $1/x$, is it as well an *approximate f-derivative* in $[1, t]$? The answer is affirmative, because the process of derivation in the scope of the IDM is very similar, reaching the equality $e^{ux} = e^{\Delta x}$, and therefore $u = \Delta x/x$, so that, it finally turns out that $f(x) = 1/x$, with $x \in [t, \infty]$, i.e., $x$ also has a lower bound specified by $t$.

---

[1] Bear in mind that the context of the IDM is discrete, so the way in which it defines the calculation of $e$ is different from usual. Thus, the presence of the limit is not necessary.

## Euler's constant

Without going into details, the rule of Barrow for the *f*-derivatives in the local calculus is given by

$$f_+(x_0) + f_+(x_0 + 1) + \cdots + f_+(x_0 + (n-1)) = F(x_0 + n) - F(x_0) \qquad \textbf{E. 3}$$

This identity is true in the interval $[1, \infty]$, obviously if the exact *f*-derivative of $F(x)$ is used. What happens if we use the approximate *f*-derivative? Replacing $1/x$ in **E. 3**, and making $x_0 = 1$, we arrive at $1 + 1/2 + 1/3 + \cdots + 1/n \approx \ln(1+n)$. Assuming that $n = t$, it is clear for the reasons stated above that the equality is not possible, since $x$ takes values in the interval $[1, t]$ and so, an *accumulative error is generated*. If we add the series up with $x \in [1, \infty]$, would the accumulated error grow indefinitely? From the previous section we know that the derivative **E. 2** is only numerically acceptable if $x \in [t, \infty]$. Therefore, as $x$ takes increasing values in the $[1, \infty]$ interval, the error generated by **E. 2** is decreasing, to become virtually non-existent once passed the threshold of $t$[1]. In cases like this, it is easy to show that the accumulative error converges on a specific value, i.e., it is a constant. Is there some easy way to calculate the accumulated error?

We know that if we apply the exact *f*-derivative in **E. 3**, with $x \in [1, \infty]$, there is no error, so a simple subtraction will provide the accumulated error. We have then that[2]

$$\gamma = \sum (f_+(x)_\sim - f_+(x)) = \sum_{x=1}^{\infty} (1/x - \ln(1 + 1/x))$$

being $\gamma = 0.577215664901...$, i.e., the famous *constant of Euler*.

What happens if we use the approximate local backward derivative *(b-derivative)* to calculate $\gamma$, instead of the *f*-derivative?

$f_-(x) = \ln(x/(x-1))$ is the exact local *b*-derivative of the function $F(x) = \ln x$. Making $f_-(x) = \ln(x/(x-1)) = u$, then $e^u = (x/(x-1))$. This equality leads to conclude that $1/e^{ux} = 1/e$, if $x \to \infty$, which again implies that $ux = 1$, so $u = 1/x$. We obtain eventually that $f_-(x)_\sim = 1/x$, and therefore the approximate derivatives (*f* and *b*) coincide, i.e., $f_+(x)_\sim = f_-(x)_\sim$.

---

[1] Since the error in each term of **E. 3** decreases as $x$ increases, the greatest contribution to the accumulative error occurs when $x = 1$, being equal to $(1 - \ln 2)$.

[2] This is the expression originally given by Euler.

Knowing the exact and approximate *b*-derivative we are able to calculate the accumulated error, which is none other than the complement to 1 of Euler's constant, i.e.,

$$\gamma = \Sigma (f_-(x) - f_-(x)_\sim) = \sum_{x=2}^{\infty}(\ln(x/(x-1)) - 1/x) = 1 - \gamma,$$

so, the accumulated error with the *b*-derivative is $\gamma = 0.422784335099...$

## The family of constants of Euler

In the wake of the results in the previous section, one wonders if there are other constants similar to $\gamma$. The answer is affirmative, because Euler's constant belongs to the family of constants associated with the functions $F(x) = (\ln x)^k$, with $k = 1, 2, 3,...$, indicated by $\gamma, \gamma_2, \gamma_3$, etc., on which we will now focus our attention.

Starting with $F(x) = (\ln x)^2$, the exact local *f*-derivative of this function is given by $f_+(x) = (\ln (x + 1))^2 - (\ln x)^2 = (\ln (x + 1) + \ln x)(\ln (x + 1) - \ln x)$, and therefore $f_+(x) = [\ln ((x + 1)x)][\ln ((x + 1)/x)]$. Let us see whether there is also an approximate local *f*-derivative.

If we do $\ln ((x + 1)x) = v$ in the exact *f*-derivative, and $\ln ((x + 1)/x) = u$, then it follows that $e^v = (x + 1)x$, so $e^v/x^2 = (1 + 1/x)$. From this equality we obtain that $(e^v/x^2)^x = e$, when $x \geq t \gg 1$, and therefore $e^{vx} = ex^{2x}$. Taking logarithms leads to $v = 1/x + 2\ln x$.

On the other hand, $\ln ((x + 1)/x) = u \Rightarrow e^u = (x + 1)/x \Rightarrow u = 1/x$, if $x \geq t \gg 1$. Taking *u* and *v* to the exact *f*-derivative, we have that

$$f_+(x)_\sim = [1/x + 2\ln x][1/x] = 1/x^2 + (2\ln x)/x,$$

that can be written as $f_+(x)_\sim = 1/x^2 + f(x)$, because $(2\ln x)/x$ is the extreme *f*-derivative of $F(x) = (\ln x)^2$, i.e., the same which provides the standard calculus.

Replacing $f_+(x)_\sim$ in the rule of Barrow (**E. 3**), if we make $x_0 = 1$ and regroup the terms of the series that arise, the result is

$$\sum_{x=1}^{n+1} \frac{1}{x^2} + 2\sum_{x=1}^{n+1} \frac{\ln x}{x} \approx (\ln(n+2))^2$$

As it happens with $k = 1$, this result is only approximate, due to the accumulative error that occurs when $x \in [1, t]$. Therefore, to obtain an accurate result we have to balance the previous expression, adding the constant $\gamma_2$ on the right side, leaving that

$$\sum_{x=1}^{n+1}\frac{1}{x^2}+2\sum_{x=1}^{n+1}\frac{\ln x}{x}=(\ln(n+2))^2+\gamma_2 \qquad \textbf{E. 4}$$

To calculate $\gamma_2$ we simply adopt the approach applied in the calculation of $\gamma$. Then, if there is an error accumulated in **E. 3** with the approximate *f*-derivative, and if there is not with the exact *f*-derivative, then we have that

$$\gamma_2 = \sum(f_+(x)_\sim - f_+(x)) = \sum_{x=1}^{\infty}\left[\left(\frac{1}{x^2}+\frac{2\ln x}{x}\right)-\ln((x+1)x)\left(\ln\left(\frac{x+1}{x}\right)\right)\right]$$

where $\gamma_2 = 1.49930237...$

Let us see the above analysis, but using the *b*-derivative.

Being much more succinct, the exact local *b*-derivative $F(x) = (\ln x)^2$ is given by $f_-(x) = \ln(x(x-1))\ln(x/(x-1))$. From it, it is not a great problem to find the approximated local *b*-derivative, which turns out to be $f_-(x)_\sim = 2\ln x/x - 1/x^2$.

The Barrow's rule for the backward derivatives is

$$f_-(x_0) + f_-(x_0+1) + \cdots + f_-(x_0+n) = F(x_0+n) - F(x_0-1) \qquad \textbf{E. 5}$$

so replacing $f_-(x)_\sim$ in this identity, and rearranging the terms after making $x_0 = 2$, we arrive at

$$2\sum_{x=2}^{n+2}\frac{\ln x}{x}-\sum_{x=2}^{n+2}\frac{1}{x^2}+\gamma_2' = (\ln(n+2))^2 \qquad \textbf{E. 6}$$

where

$$\gamma_2' = \sum(f_-(x)-f_-(x)_\sim) = \sum_{x=1}^{\infty}\left[\ln(x(x-1))\left(\ln\left(\frac{x}{x-1}\right)\right)-\left(\frac{2\ln x}{x}-\frac{1}{x^2}\right)\right]$$

On the other hand, subtracting **E. 6** from **E. 4**, comes to

$$1+2\sum_{x=2}^{n+1}\frac{1}{x^2}+\frac{1}{(n+2)^2}-\frac{2\ln(n+2)}{n+2} = \gamma_2+\gamma_2'$$

Now, if $n \to \infty$, or also, if we say that $x \in [1, \infty]$, then we can write[1]

$$\sum_{x=2}^{\infty}\frac{1}{x^2} = \frac{\gamma_2+\gamma_2'-1}{2}$$

Adding 1 to both sides, we have that

---

[1] Outside the scope of the IDM these deductions may not seem very orthodox, but in the discrete context of the IDM they are completely rigorous.

$$\sum_{x=1}^{\infty}\frac{1}{x^2} = \frac{\gamma_2 + \gamma_2' + 1}{2} = \frac{\pi^2}{6} \Rightarrow \gamma_2 + \gamma_2' = \frac{\pi^2}{3} - 1$$

To conclude, we will now calculate $\gamma_3$ from the function $F(x) = (\ln x)^3$, but with a more abbreviated development, because it is similar to the previous ones.

We start looking for the exact local *f*-derivative in the $[1, \infty]$ interval, that has the general approach $f_+(x) = (\ln (x + 1))^3 - (\ln x)^3 = (\ln (x + 1) - \ln x)( (\ln (x + 1))^2 - \ln (x + 1)\ln x + (\ln x)^2 )$, and so,

$$f_+(x) = \ln ((x + 1)/x)( (\ln (x + 1))^2 + \ln (x + 1)\ln x + (\ln x)^2 )$$

If $\ln ((x + 1)/x) = u$, $\ln (x + 1) = t$, and $\ln (x + 1)\ln x = v$, we have that $u = 1/x$, $t^2 = (1/x + \ln x)^2$, and $v = (\ln x)/x + (\ln x)^2$. Taking *u*, *t*, y *v* to the $f_+(x)$ definition, it turns out that $f_+(x)_\sim = 1/x^3 + 3(\ln x)/x^2 + 3(\ln x)^2/x$.

This is the approximate local *f*-derivative of $F(x) = (\ln x)^3$, which can also be written as $f_+(x)_\sim = 1/x^3 + 3(\ln x)/x^2 + f(x)$, where $f(x)$ is the extreme derivative of $F(x) = (\ln x)^3$. Substituting $f_+(x)_\sim$ in **E. 3**, we have that

$$\sum_{x=1}^{n}\frac{1}{x^3} + 3\sum_{x=1}^{n}\frac{\ln x}{x^2} + 3\sum_{x=1}^{n}\frac{(\ln x)^2}{x} \approx (\ln x)^3$$

The accumulated error on this occasion is

$$\gamma_3 = \sum(f_+(x)_\sim - f_+(x)) = 3.9856304...,$$

and so

$$\sum_{x=1}^{n}\frac{1}{x^3} + 3\sum_{x=1}^{n}\frac{\ln x}{x^2} + 3\sum_{x=1}^{n}\frac{(\ln x)^2}{x} = (\ln x)^3 + \gamma_3 \quad \text{E. 7}$$

The approach of the *b*-derivative in this case is $f_-(x) = (\ln x)^3 - (\ln (x - 1))^3$, and from it, repeating the steps given above, it turns out that $f_-(x) = \ln (x/(x - 1))((\ln x)^2 + \ln x \ln (x - 1) + (\ln (x - 1))^2)$.

Now, if $\ln (x/(x - 1)) = u$, $\ln (x - 1) = t$, and $\ln x \ln (x - 1) = v$ we have that $u = 1/x$, $t^2 = (\ln x - 1/x)^2$, and $v = (\ln x)^2 - (\ln x)/x$. Replacing these three equalities in the previous one we obtain

$$f_-(x)_\sim = 3(\ln x)^2/x - 3(\ln x)/x^2 + 1/x^3$$

Taking this *b*-derivative to **E. 5**, it comes to

$$3\sum_{x=2}^{n+2}\frac{(\ln x)^2}{x} - 3\sum_{x=2}^{n+2}\frac{\ln x}{x^2} + \sum_{x=2}^{n+2}\frac{1}{x^3} + \gamma_3' = (\ln (n + 2))^3$$

where $\gamma_3' = \sum (f_+(x) - f_+(x)_-) = 2.6396589...$

Finally, subtracting this expression from **E. 7**, we obtain

$$\gamma_3 + \gamma_3' = 6\sum_{x=1}^{n} \frac{\ln x}{x^2} + 1$$

and replacing the summatory by the definite integral, we have finally that $\gamma_3 + \gamma_3' = 7 - 2\lambda_1$, being $\lambda_1$ a meaningful constant related to the terms of **E. 7**, whose value is 0.1873553...

We conclude here the analysis of the family of γ. We have seen that the calculation of Euler's constants thickens gradually, as the exponent $k$ of the function $F(x) = (\ln x)^k$ grows. However, everything seems to indicate that it will always be possible to find an approximate derivative of $F(x)$, regardless of the value of $k$, which would mean that the total membership in the family of constants of Euler would be unlimited.

## Conclusions

Although the use of local derivatives in the study of the family of the constants of Euler is not enough to evaluate the possibilities of the IDC thoroughly, this particular analysis illustrates quite well the characteristics and capabilities of the isodimensional calculus. We have seen how the choice of the level of observation in the range [$S_0$, $S_m$] determines the value of the variation $\Delta x$, which in turn sets the characteristics of the derivative functions and integrals (not seen), that is to say, either they resemble the expressions that we find in the finite differences methods (when $m = 0$), or they are equal to those provided by the standard calculus, when $m \to \infty$. In any case, regardless of the level of valuation $S_0$ we adopt, the concepts and techniques of integration and derivation are more or less the same. Proof of this is that the quotient of increments **E. 1** is always the same, whatever the chosen observation level, because the discrete context provided by the IDM does not require the presence of the limit in any case, something that, on the other hand, facilitates the algebraic and deductive processes.